\documentclass[12pt]{amsart}
\usepackage{verbatim}
\usepackage{amsmath}
\usepackage{amsthm}
\usepackage{amssymb}
\usepackage{enumerate}

\newif\ifdeveloping

\newtheorem{theorem}{Theorem}
\newtheorem{proposition}[theorem]{Proposition}

\newtheorem{corollary}[theorem]{Corollary}

\newtheorem{problem}[theorem]{Problem}

\theoremstyle{definition}

\newtheorem{definition}[theorem]{Definition}

\theoremstyle{remark}
  
{}

\ifdeveloping
\usepackage[notref,notcite]{showkeys}
\fi

\newcommand{\prtime}{{\count0=\time\divide\count0 by 60
\count1=-\count0\multiply\count1 by 60 \advance\count1 by \time
\the\count0:\the\count1} }

\def\myheads#1;#2;{
\pagestyle{myheadings} \markboth{{\sc\hfill
#1\hfill\protect\makebox[0cm][r]{\rm\today; \prtime}}}
{{\sc\protect\makebox[0cm][l]{\rm\today;\ \prtime}\hfill #2\hfill}}
\thispagestyle{myheadings} }

\newcommand{\subs}{\subset}

\newcommand{\dom}{\operatorname{dom}}

\def\<{\left\langle}
\def\>{\right\rangle}
\def\cf{\operatorname{cf}}
\def\cof{\operatorname{cof}}

\def\br#1;#2;{\bigl[ {#1} \bigr]^ {#2} }
\def\bc#1;#2;{\bigl( {#1} \bigr)^ {#2} }

\theoremstyle{plain}

\begin{document}
\author[I. Juh\'asz]{Istv\'an Juh\'asz}
\address{Alfr{\'e}d R{\'e}nyi Institute of Mathematics}
\email{juhasz@renyi.hu}

\author[L. Soukup]{Lajos Soukup }
\address{Alfr{\'e}d R{\'e}nyi Institute of Mathematics }
\email{soukup@renyi.hu}

\author[Z. Szentmikl\'ossy]{Zolt\'an Szentmikl\'ossy}
\address{E\"otv\"os Lor\'ant University, Department of  Analysis,
1117 Budapest,
 P\'azm\'any P\'eter s\'et\'any 1/A,
 Hungary}
\email{zoli@renyi.hu}
\thanks{Research supported
by OTKA grant no.\ 61600.}



\title[Point-countable $\pi$-bases]{First countable spaces without point-countable $\pi$-bases}

\begin{abstract}
We answer several questions of V. Tka\v cuk from [{\em
Point-countable $\pi$-bases in first countable and similar spaces},
Fund. Math. 186 (2005), pp. 55--69.] by showing that
\begin{enumerate}
\item there is a ZFC example of a first countable, 0-dimensional Hausdorff space
with no point-countable $\pi$-base (in fact, the order of any
$\pi$-base of the space is at least $\aleph_\omega$);
\smallskip
\item if there is a $\kappa$-Suslin line then there is a first
countable GO space of cardinality $\kappa^+$ in which the order of
any $\pi$-base is at least $\kappa$;
\smallskip
\item it is consistent to have a first
countable, hereditarily Lindel\" of regular space having uncountable
$\pi$-weight and $\omega_1$ as a caliber (of course, such a space
cannot have a point-countable $\pi$-base).
\end{enumerate}
\end{abstract}

\maketitle

\section{Introduction}

V. Tka\v cuk in \cite{T} has recently proved under CH that any first
countable Hausdorff space that is Lindel\" of or CCC has a
point-countable $\pi$-base. (Actually, in \cite{T} all spaces are
assumed to be Tychonov, but the proof only needs Hausdorff.) Tka\v
cuk's motivation was to extend (at least partially) \v Sapirovskii's
celebrated ZFC result saying that any countably tight compactum has
a point-countable $\pi$-base, from compact spaces to Lindel\" of
ones. So it was natural to ask if his use of CH was necessary. To
Tka\v cuk's amazement (and embarrassment), he could come up with no
example of a first countable space without a point-countable
$\pi$-base, not even consistently.

Our aim here is to remedy this embarrassing situation and provide
not only ZFC (and several consistent)  examples of first countable
(Tychonov) spaces without point-countable $\pi$-bases but also
examples which show that, indeed, Tka\v cuk's CH results cannot be
proved in ZFC alone. In this manner we succeeded in answering 7 of
the 12 questions that were listed at the end of \cite{T}.

In what follows, we shall use the notation and terminology of
\cite{J}. In particular, $\pi sw(X)$ denotes the {\em
$\pi$-separating weight} of $X$, that is the minimum order of a
$\pi$-base of the space $X$, see p. 74 of \cite{J}. Note that $\pi
sw(X) \le \omega$ is then equivalent to the statement: $X$ has a
point-countable $\pi$-base.

\section{ZFC examples}

The key to Tka\v cuk's above mentioned CH results in \cite{T} was
his theorem 3.1 which says that if $X$ has countable tightness and
$\pi$-character, moreover $d(X) \le \omega_1$ then $\pi sw(X) \le
\omega$. In his list of problems (problem 4.11), Tka\v cuk himself
asked if the assumption of countable tightness could be omitted
here. It is immediate from our next result that this question has an
affirmative answer.

\begin{theorem}\label{pich}
Let $X$ be any topological space with $d(X) \le \pi \chi(X)^+$. Then
$\,\,\pi sw(X) \le \pi \chi(X)$.
\end{theorem}

\begin{proof}
Let us set $\pi \chi(X) = \kappa$. If $d(X) \le \kappa$ then we even
have $\pi(X) = \kappa$. So we may assume $d(X) = \kappa^+$ and, as
is well-known, we may then fix a dense set $D =\{ x_\alpha : \alpha
< \kappa^+ \}$ that is left-separated in this well-ordering. This
means that for every $\alpha < \kappa^+$ there is a neighbourhood
$U_\alpha$ of $x_\alpha$ with $$\{ x_\beta : \beta < \alpha \} \cap
U_\alpha = \emptyset.$$ Let us now fix a local $\pi$-base
$\mathcal{B}_\alpha$ of the point $x_\alpha$ such that
$|\mathcal{B}_\alpha| \le \kappa$ and $B \subset U_\alpha$ whenever
$B \in \mathcal{B}_\alpha$. Then $\mathcal{B} = \bigcup \{
\mathcal{B}_\alpha : \alpha < \kappa^+ \}$ is clearly a $\pi$-base
of $X$ such that for every $x_\beta \in D$ we have $$ord(x_\beta,
\mathcal{B}) = |\{ B \in \mathcal{B} : x_\beta \in B \}| \le
\kappa.$$

We claim that then we have $$ord(\mathcal{B}) = \sup \{
ord(x,\mathcal{B}): x \in X \} \le \kappa$$ as well. Assume, on the
contrary, that $\,ord(x,\mathcal{B}) = \kappa^+$ for some point $x
\in X.$ Since $\pi\chi(x,X) \le \kappa$, this implies that there are
$\kappa^+$ many members of $\mathcal{B}$ (containing $x$) that
include a fixed non-empty open set $V$. This, however, is impossible
because $D \cap V \ne \emptyset.$
\end{proof}

We may now turn to our main aim that is to produce, in ZFC,  first
countable spaces without point-countable $\pi$-bases.

\begin{theorem}\label{ZFC}
There is a first countable, 0-dimensional Hausdorff (hence Tychonov)
space $X$ with $\pi sw(X) \ge \aleph_\omega$.
\end{theorem}

\begin{proof}
The underlying set of our space is  $X = \prod \{ \omega_n : n <
\omega \}$. For $f,g \in X$ we write $f \le g$ to denote that $f(n)
\le g(n)$ for all $n < \omega.$ The topology $\tau$ that we shall
consider on $X$ will be generated by all sets of the form $U_n(f)$
(with $f \in X$ and $n < \omega$), where $$U_n(f) = \{ g \in X : f
\le g \mbox{ and } f\upharpoonright n = g\upharpoonright n \}.$$
Note that if $g \in U_n(f)$ then $U_n(g) \subs U_n(f)$, and if $g
\notin U_n(f)$ then there is a $k < \omega$ such that $U_k(g) \cap
U_n(f) = \emptyset$. It follows immediately from this that, for any
$f \in X$, the family $\{ U_n(f) : n < \omega\}$ forms a clopen
neighbourhood base of $f$ with respect to the topology $\tau$,
consequently the space $\langle X,\tau \rangle$ is indeed first
countable, 0-dimensional, and Hausdorff.

It is also easy to see from the definitions that if
$\{U_{n_{\alpha}}(f_{\alpha}) : \alpha < \kappa\}$ is a $\pi$-base
of $\tau$ then $\{f_{\alpha} : \alpha < \kappa\}$ must be cofinal in
the partial order $\langle X, \le \rangle$. But it is well-known
that the latter has cofinality greater than $\aleph_\omega$,
consequently we have $\pi(X) > \aleph_\omega.$ (Actually, it is easy
to see that $\pi(X) = \cf\big(\langle X, \le \rangle \big)$ but we
shall not need this.)

Next we claim that for any $k < \omega$ the pair
$(\aleph_{\omega+1}, \aleph_k)$ is a {\em pair caliber} of the space
$X$, i. e. among any $\aleph_{\omega+1}$ open sets one can find
$\aleph_k$ whose intersection is non-empty. Indeed, we may assume
without any loss of generality that we have a family of basic open
sets of the form $\{U_n(f) : f \in F\}$ where $F \in
[X]^{\aleph_{\omega+1}}$ and $n > k$ is fixed. We may also assume
that $f \upharpoonright n = \sigma$ for a fixed $\sigma \in \prod_{i
< n}\omega_i$ whenever $f \in F$. Now let $G \subs F$ with $|G| =
\aleph_k$ then there is $g \in X$ with $g \upharpoonright n =
\sigma$, moreover $f(i) < g(i)$ for all $i \ge n$ and $f \in G$. But
then we clearly have $g \in \bigcap \{ U_n(f) : f \in G \}$.

Now, putting together the previous two paragraphs we can immediately
conclude that the order of any $\pi$-base of $\langle X,\tau
\rangle$, being of cardinality  $> \aleph_\omega$, must be at least
$\aleph_\omega,\,\,$ i. e. $\pi sw(X) \ge \aleph_\omega$.

\end{proof}

The cardinality of our above example is $\aleph_\omega^{\aleph_0}$
that is much larger than the optimal value $\aleph_2$ permitted by
theorem \ref{pich}. So it is natural to raise the question if we
could find other examples of smaller cardinality. As it turns out,
we can do slightly better by choosing an appropriate subspace $Y$ of
the space $X$ from theorem \ref{ZFC}. First, however, we need to fix
some notation. For $f, g \in X = \prod \{ \omega_n : n < \omega \}$
we write $f <^* g$ to denote that $|\{ n < \omega : f(n) \ge g(n)
\}|$ is finite, i. e. $f$ is below $g$ modulo finite. Similarly, we
write $f =^* g$ to denote that $|\{ n < \omega : f(n) \ne g(n) \}|$
is finite. Finally, we recall that there is in $X$ a transfinite
sequence of order type $\omega_{\omega+1}$ that is increasing with
respect to $<^*$.

\begin{theorem}\label{subZFC}
Let $\,\{ f_\alpha : \alpha < \omega_{\omega+1} \} \subs X$ be an
increasing sequence with respect to $<^*$ and set $$Y = \{ f \in X
:\, \exists\,\alpha < \omega_{\omega+1} \mbox{ with } f =^* f_\alpha
\}.$$ Then the subspace $Y$ of $X$, with the subspace topology
inherited from $\tau$, also satisfies $\pi sw(Y) \ge \aleph_\omega.$
\end{theorem}

\begin{proof}
The proof is very similar to that of theorem \ref{ZFC}. First we
note that, trivially, again we have $\pi(Y) > \aleph_\omega.$ Next,
we show that $(\aleph_{\omega+1}, \aleph_k)$ is a  pair caliber of
$Y$ for each $k < \omega.$ To see this, we let again consider any
family $\{U_n(f) : f \in F\}$ where $F \in [Y]^{\aleph_{\omega+1}}$
and $n > k > 0$, moreover $f \upharpoonright n = \sigma$ for a fixed
$\sigma \in \prod_{i < n}\omega_i$ whenever $f \in F$. Let us choose
any subset $G \subs F$ with $|G| = \aleph_k$, then there is an
ordinal $\alpha < \omega_{\omega+1}$ such that $g <^* f_\alpha$ for
all $g \in G$. Clearly, then we may find an integer $m \ge n$ such
that the set $$G^* = \{ g \in G : \,\forall\,i \ge m\,\,\big(g(i) <
f_\alpha(i) \big) \}$$ also has cardinality $\aleph_k$.

Note that if $n \le j < m$ then $\{g(j) : g \in G^*\}$ is bounded in
$\omega_j$, hence we may find a function $f \in Y$ such that $f
\upharpoonright n = \sigma$, if $n \le j < m$ then $g(j) < f(j)$ for
all $g \in G^*$, moreover $f(i) = f_\alpha(i)$ whenever $m \le i <
\omega$. But then we have $f \in \bigcap \{ U_n(g) : g \in G^* \}
\cap Y$.
\end{proof}

We were unable to produce a ZFC example of a first countable space
without a point-countable $\pi$-base of cardinality less than
$\aleph_{\omega+1}$. This leads us to the following intriguing open
question.

\begin{problem}\label{kis}
Is there, in ZFC, a first countable (Tychonov) space of cardinality
less than $\aleph_\omega$ that has no point-countable $\pi$-base?
\end{problem}

Actually, at this point we do not even have such an example of
cardinality $\aleph_\omega$. We conjecture, however, that having
such an example is equivalent to having one of size $<
\aleph_\omega$. In fact, we could verify this conjecture under the
assumption $2^{\aleph_1} < \aleph_\omega$.

\begin{theorem}\label{ekv}
Assume that $2^{\aleph_1} < \aleph_\omega$ and $X$  is a first
countable space of cardinality $\aleph_\omega$. If every subspace of
$X$ of cardinality $< \aleph_\omega$ has a point-countable
$\pi$-base then so does $X$.
\end{theorem}

\begin{proof}
Let us start by giving a (very natural) definition. A family
$\mathcal{B}$ of non-empty open sets in $X$ is said to be an {\em
outer $\pi$-base} of a subspace $Y \subs X$ if for every open set
$U$ with $U \cap Y \ne \emptyset$ there is a member $B \in
\mathcal{B}$ such that $B \subs U$. We claim that, under the
assumptions of our theorem, every subspace of $X$ of cardinality $<
\aleph_\omega$ actually has a point-countable outer $\pi$-base. Thus
if we have $X = \bigcup_{n < \omega} Y_n$ where $|Y_n| <
\aleph_\omega$ for all $n < \omega$ and $\mathcal{B}_n$ is a
point-countable outer $\pi$-base of $Y_n$ in $X$ then $\bigcup_{n <
\omega} \mathcal{B}_n$ is clearly a point-countable $\pi$-base of
$X$.

To prove the above claim let us consider an $\omega_1$-closed
elementary submodel $M$ with $|M| < \aleph_\omega$ of a "universe"
$H(\theta)$. (As usual, here $H(\theta)$ is the collection of all
sets of hereditary cardinality $< \theta$, and  $M$ is
$\omega_1$-closed means that $[M]^{\le \omega_1} \subs M$.) The
regular cardinal $\theta$ is chosen so large that $H(\theta)$ (and
also $M$) contains $X$ and everything else that is relevant, e. g. a
map $\mathcal{V}$ that assigns to every point $x \in X$ a countable
open neighbourhood base $\mathcal{V}_x$. Clearly, $2^{\aleph_1} <
\aleph_\omega$ implies that for every $Y \in [X]^{< \aleph_\omega}$
there is such an elementary submodel $M$ with $Y \subs M$.
Consequently, our claim will be proven if we show that $X \cap M$
has a point-countable outer $\pi$-base in $X$ whenever $M$ is like
above.

To see this, note first that for every point $x \in X \cap M$ we
have $\mathcal{V}_x \in M$ and hence $\mathcal{V}_x \subs M$ as
well. Consequently $\mathcal{V}_M = \cup \{ \mathcal{V}_x : x \in X
\cap M \} \subs M$ is an outer base of $X \cap M$ in $X$, hence we
may choose a subfamily $\mathcal{B} \subs \mathcal{V}_M$ such that
$\mathcal{B} \upharpoonright M = \{ B \cap M : B \in \mathcal{B} \}$
is a point-countable $\pi$-base of the subspace $X \cap M$.

It suffices to show now that $\mathcal{B}$ is a point-countable
outer $\pi$-base of $X \cap M$ in $X$. Indeed, $\mathcal{B}$ is
point-countable for if $\mathcal{U} \in [\mathcal{B}]^{\omega_1}$
then we have $\mathcal{U} \in M$ because $M$ is $\omega_1$-closed,
and thus $\cap \mathcal{U} \ne \emptyset$ would imply $\cap
\mathcal{U} \cap M \ne \emptyset$, contradicting that $\mathcal{B}
\upharpoonright M$ is point-countable. (Here we also used the fact
that, by elementarity, the correspondance $B \mapsto B \cap M$ is
one-to-one on $\mathcal{B} \subs M$.) Finally, $\mathcal{B}$ is an
outer $\pi$-base of $X \cap M$ because if $U$ is open with $x \in U
\cap M \ne \emptyset$ then there is $V \in \mathcal{V}_x \subs M$
with $V \subs U$, hence if $B \in \mathcal{B}$ with $B \cap M \subs
V \in M$ then we also have $B \subs V \subs U$.
\end{proof}

\section{Examples from higher Suslin lines}

We start this section by giving a theorem that, not surprisingly,
will turn out to be very useful in finding (first countable) spaces
without point-countable $\pi$-bases.

\begin{theorem}\label{dis}
Assume that $X$ is a topological space which has a $\pi$-base
$\mathcal{B}$ such that $ord(\mathcal{B})^+ < d(X)$. Then $X$ has a
discrete subspace $D$ with $|D| \ge d(X).$
\end{theorem}

\begin{proof}
Let us first choose a point $x_B \in B$ from each $B \in
\mathcal{B}$. Clearly, then the set $S = \{ x_B : B \in \mathcal{B}
\}$ is dense in $X$, hence we have $$|S| \ge d(X) >
ord(\mathcal{B})^+.$$ We now define a set mapping $F$ on $S$ by the
following stipulation: For any point $x \in S$ let us put $$F(x) =
\{ x_B \in S : x \in B \} \in [S]^{\le ord(\mathcal{B})}.$$ By
Hajnal's set mapping theorem (see \cite{H}) then there is a free set
$D \subs S$ for the set mapping $F$ with $|D| = |S|$. This means
that for every $x \in D$ we have $D \cap F(x) \subs \{x\}.$ But
every member of $D$ is of the form $x_B$ for some $B \in
\mathcal{B}$, and we claim that for this point we have $B \cap D =
\{x_B\}.$ Indeed, $x_B \in B \cap D$ is obvious, and if $x \in D$ is
different from $x_B$ then $x_B \notin F(x)$ implies $x \notin B$.
Consequently, $D$ is as required.
\end{proof}

It is an immediate consequence of theorem \ref{dis} that a space $X$
satisfying $s(X) < d(X) \ge \omega_2$ cannot have a point-countable
$\pi$-base. Unfortunately, we do not know if there is in ZFC a first
countable Tychonov space like that. (Recall that the solution of
Tka\v cuk's problems from require Tychonov examples.)

If, however, we are satisfied with Hausdorff examples then we are
much better off.  In fact it was shown in \cite{HJ} that there is a
natural left-separated refinement $\sigma$ of the euclidean topology
$\tau$ on the real line $\mathbb{R}$ that is first countable and
hereditarily Lindel\"of, moreover the density of $\sigma$ is equal
to $\cf(\mathfrak{c})$. Consequently, if $\cf(\mathfrak{c}) >
\omega_1$ then, by theorem \ref{dis}, the space $\langle \mathbb{R},
\sigma \rangle$ has no point-countable $\pi$-base. This is of
interest because it shows that, at least for Hausdorff spaces, Tka\v
cuk's CH result that was mentioned in the introduction is not valid
without CH, for $\langle \mathbb{R}, \sigma \rangle$ is both
Lindel\"of and CCC. Actually, it is very easy to show that, assuming
$\cf(\mathfrak{c}) > \omega_1$ again, something stronger than CCC
can be established, namely that $\omega_1$ is a caliber of this
space and, for Hausdorff spaces, this settles one more question of
Tka\v cuk from \cite{T}. In the next section we shall produce
(consistent) Tychonov examples with these properties but that will
require more work.

Next we shall consider higher Suslin lines, these are ordered spaces
whose spread is less than their density. More precisely, we shall
consider first countable variations of them that retain this
property. For different purposes, this construction had been already
used in theorem 1.1 of \cite{JHB}, although there CH was
additionally assumed.

Let $\kappa$ be an infinite cardinal.  We shall call a continuous
linear order $\langle L, < \rangle$, equipped with the order
topology generated by $<\,$, a $\kappa$-Suslin line if there are no
more than $\kappa$ disjoint open intervals in $L$ (i. e. $c(L) \le
\kappa$), although the density $d(L)$ of $L$ is bigger than
$\kappa$. (It is known that the existence of a $\kappa$-Suslin line
is equivalent to the existence of a $\kappa$-Suslin tree, but this
will be irrelevant for us.) Thus, an ordinary Suslin line is the
same as an $\omega$-Suslin-line and by a higher Suslin line we mean
a $\kappa$-Suslin line where $\kappa > \omega$.

The main result of this section is the following theorem that, in
particular, yields us a consistent example of a first countable
GO-space without a point-countable $\pi$-base of the minimum
possible cardinality $\omega_2$. (Recall that GO-spaces, or
generalized ordered spaces, are the subspaces of linearly ordered
spaces.)

\begin{theorem}\label{sus}
If there is a $\kappa$-Suslin line $\langle L, < \rangle$ then there
is a first countable GO-space $X$ with $|X| = \kappa^+$ and $\pi
sw(X) = \kappa$.
\end{theorem}

\begin{proof}
Let $Z$ be the set of all those points $x \in L$ that have
left-character $\omega$, that is the open half line $(\leftarrow,x)$
has cofinality $\omega$ with respect to $<$. Since $\langle L, <
\rangle$ is continuous, $Z$ is clearly dense in $L$. It follows that
$d(Z) = d(L) = \kappa^+$ because $d(L) \le c(L)^+$ holds for any
linearly ordered space $L$, see e. g. \cite{B}. Now let $X$ be any
dense subspace of $Z$ (and hence of $L$) with $|X| = \kappa^+$.

We consider $X$ with the left-Sorgenfrey topology $\sigma$, i. e.
for any $x \in X$ the half-open intervals $(y,x]$ form a
$\sigma$-local base. Then $\sigma$ is clearly finer than the order
topology on $X$, hence the density of $\langle X,\sigma \rangle$
must be bigger than $\kappa$. It is clear from the definition that
$\sigma$ is a first countable topology.

Also, $\langle X,\sigma \rangle$ is a GO space because it is
homeomorphic to the subspace topology on $X \times \{ 0 \}$
inherited from the order topology on $L \times 2$ taken with the
lexicographic order. Finally, we clearly have $c(X,\sigma) = c(L)$,
moreover $s(X,\sigma) = c(X,\sigma)$ is known to hold for GO-spaces,
see e. g. 2.23 of \cite{J}. Consequently we have $s(X) \le \kappa <
d(X)$ and so theorem \ref{dis} implies $\pi sw(X) \ge \kappa$. But
by theorem \ref{pich} we must have then $\pi sw(X) = \kappa$.
\end{proof}
In particular, the existence of an $\omega_1$-Suslin line implies
that of a first countable GO space of cardinality $\omega_2$ without
a point-countable $\pi$-base.

Finally, we mention here the curious fact that it is an outstanding
open question of set theory whether one can find a model of ZFC that
does not contain any higher Suslin line. Consequently there is a
chance that theorem \ref{sus} yields us a ZFC example of a first
countable GO space with no point-countable $\pi$-base.

\section{Examples from subfamilies of $\mathcal{P}(\omega)$ }

In this section we are going to introduce a (quite simple but
apparently new) way of constructing first countable, 0-dimensional
Hausdorff topologies on subfamilies of $\mathcal{P}(\omega)$, the
power set of $\omega$. Then we shall use some of the spaces obtained
in this manner to present examples that demonstrate the necessity of
the use of CH in Tka\v cenko's results mentioned in the
introduction.

We start with fixing some notation and terminology. We shall use
$\mathcal{U}$ to denote the family of all co-finite subsets of
$\omega$. For a given family $\mathcal{I} \subs \mathcal{P}(\omega)$
and for $I \in \mathcal{I}$ and $U \in \mathcal{U}$ we put
$$[I,U)_{\mathcal{I}} = \{ J \in \mathcal{I} : I \subs J \subs
U\}.$$ If $\mathcal{I} = \mathcal{P}(\omega)$ then we shall omit the
subscript.

Finally, we say that the family $\mathcal{I} \subs
\mathcal{P}(\omega)$ is {\em stable} if $I \in \mathcal{I}$ and $I
=^* J$ for $J \subs \omega$ imply $J \in \mathcal{I}$ as well. (Of
course, here $I =^* J $ means that $I$ and $J$ are equal mod finite,
i. e. their symmetric difference $I \Delta J$ is finite.)

\begin{definition}
Let us fix a family $\mathcal{I} \subs \mathcal{P}(\omega)$. We
shall denote by $\tau_{\mathcal{I}}$ the topology on $\mathcal{I}$
generated by all sets of the form $[I,U)_{\mathcal{I}}\,$, where $I
\in \mathcal{I}$ and $U \in \mathcal{U}\,$, and by $X_{\mathcal{I}}$
the space $\langle \mathcal{I}, \tau_{\mathcal{I}} \rangle$.
\end{definition}

Of course, $X_{\mathcal{I}}$ is identical with the appropriate
subspace of the maximal such space $X_{\mathcal{P}(\omega)}$. A few
basic (pleasant) properties of the spaces $X_{\mathcal{I}}$ are
given by the following proposition.

\begin{proposition}
The spaces $X_{\mathcal{I}}$ are first countable, 0-dimensional and
Hausdorff.

\end{proposition}

\begin{proof}
Clearly, it suffices to show this for $\mathcal{I} =
\mathcal{P}(\omega)$ because all three properties are inherited by
subspaces.

Now, observe first that if $J \in [I,U) \cap [I',U')$ then $$J \in
[J,U \cap U') \subs [I,U) \cap [I',U'),$$ hence the ``intervals"
$[I,U)$ form an open basis of $\tau_{\mathcal{P}(\omega)}$, moreover
$\{ [I,U) : I \subs U \in \mathcal{U}\}$ forms a countable
neighbourhood base of the point $I$ of $X_\mathcal{I}$.

Next, if $J \notin [I,U)$ then either $J \backslash U \ne \emptyset$
and then $[J,\omega) \cap [I,U) = \emptyset$, or $J \subs U$ and $I
\backslash J \ne \emptyset$. But in the latter case we may pick an
$n \in I \backslash J$ and then we have $J \subs U \backslash
\{n\}$, moreover $[J,U \backslash \{n\}) \cap [I,U) = \emptyset$
because $n \in I$. This means that all basic open sets $[I,U)$ are
also closed, that is $X_{\mathcal{P}(\omega)}$ is indeed
0-dimensional.

Finally, for every $I \in \mathcal{P}(\omega)$ we clearly have
$$\bigcap \{ [I,U) : I \subs U \in \mathcal{U} \} = \{ I \},$$
implying that $X_\mathcal{P}(\omega)$ is also Hausdorff.
\end{proof}

For any family $\mathcal{I} \subs \mathcal{P}(\omega)$ we shall
denote by $\cof(\mathcal{I})$ the cofinality of the partial order
$\< \mathcal{I}, \subs \>$. Also, we say that a cardinal number
$\kappa$ is a {\em set caliber} of $\mathcal{I}$ if for every
subfamily $\mathcal{J} \in [\mathcal{I}]^{\kappa}$ there are
$\mathcal{K} \in [\mathcal{J}]^{\kappa}$ and $I \in \mathcal{I}$
such that $\cup \mathcal{K} \subs I$ or, less formally, among any
$\kappa$-many members of $\mathcal{I}$ there are $\kappa$-many that
have an upper bound in $\mathcal{I}$. We now connect these concepts
concerning $\mathcal{I}$ with properties of the associated space
$X_{\mathcal{I}}$.

\begin{proposition}\label{ekv}
For any subfamily $\mathcal{I} \subs \mathcal{P}(\omega)$ we have
\begin{enumerate}[(i)]
\item $d(X_{\mathcal{I}}) = \cof(\mathcal{I})\cdot \omega\,$;

\item if $\mathcal{I}$ is stable and $\kappa$ is a cardinal with
$\cf(\kappa) > \omega$ then $\kappa$ is a caliber of the space
$X_{\mathcal{I}}$ if and only if $\kappa$ is a set caliber of the
family $\mathcal{I}$.
\end{enumerate}
\end{proposition}

\begin{proof}
The proof of (i) and the left-to-right direction of (ii) follows
immediately from the fact that $\mathcal{K} \subs \mathcal{I}$ has
an upper bound in $\mathcal{I}\,$ iff $\,\bigcap \{
[I,\omega)_{\mathcal{I}} : I \in \mathcal{K} \} \ne \emptyset\,.$ To
see the other direction, assume that $\kappa$ is a set caliber of
the family $\mathcal{I}$ and consider a family $\mathcal{B}$ of
$\kappa$-many basic open sets. Since $\cf(\kappa) > \omega$ we may
assume that  $\mathcal{B} = \{ [I,U) : I \in \mathcal{J} \}$ for
$\mathcal{J} \in [\mathcal{I}]^{\kappa}$ and a fixed $U \in
\mathcal{U}$. By our assumption there is a $\mathcal{K} \in
[\mathcal{J}]^{\kappa}$ which has an upper bound $K \in
\mathcal{I}$. Clearly, then $K \cap U \in \mathcal{I}$ as
$\mathcal{I}$ is stable and consequently $$K \cap U \in \bigcap
\big\{ [I,U) : I \in \mathcal{K} \big\}.$$

\end{proof}

After these preparatory propositions we can now present a result
that will allow us to obtain nice examples of first countable spaces
without point-countable $\pi$-bases.

\begin{theorem}\label{ter}
Assume that $\mathcal{I} \subs \mathcal{P}(\omega)$ is stable,
$\cof(\mathcal{I}) > \omega$, and $\omega_1$ is a set caliber of
$\mathcal{I}$. Then $\pi sw(X_{\mathcal{I}}) > \omega$.
\end{theorem}

\begin{proof}
Since $X_{\mathcal{I}}$ is first countable, and by (i) of
proposition \ref{ekv}, we have $$\pi(X_{\mathcal{I}}) =
d(X_{\mathcal{I}}) = \cof(\mathcal{I}) > \omega.$$ But, in view of
part (ii) of proposition \ref{ekv}, $\omega_1$ is a caliber of
$X_{\mathcal{I}}$, consequently no $\pi$-base of $X_{\mathcal{I}}\,$
can be point-countable.
\end{proof}

\begin{corollary}\label{inc}
Assume that there is a mod finite strictly increasing
$\omega_2$-sequence in $\mathcal{P}(\omega)$. Then there is a first
countable, 0-dimensional and Hausdorff space of cardinality
$\omega_2$ and having $\omega_1$ as a caliber. In particular,
$MA_{\omega_1}$ implies the existence of such a space.
\end{corollary}

\begin{proof}
Let $\{ A_\alpha : \alpha < \omega_2 \} \subs \mathcal{P}(\omega)$
be a mod finite strictly increasing $\omega_2$-sequence, i. e. we
have $|A_\alpha \backslash A_\beta| < \omega$ and $|A_\beta
\backslash A_\alpha| = \omega$ whenever $\alpha < \beta < \omega_2$.
It is obvious that the family $$\mathcal{I} = \{ I \subs \omega :
\exists\,\alpha < \omega_2 \mbox{ with } I =^* A_\alpha  \}$$ is
stable and satisfies $|\mathcal{I}| = \cof(\mathcal{I}) = \omega_2$.
Next, we claim that $\mathcal{I}$ has $\omega_1$ as a set caliber.

To see this, consider any family $\mathcal{J} = \{ I_\alpha : \alpha
\in a \} \subs \mathcal{I}$ where $a \in [\omega_2]^{\omega_1}$ and
$I_\alpha =^* A_\alpha$ for all $\alpha \in a$ and pick $\beta <
\omega_2$ such that $a \subs \beta$. Then $|A_\alpha \backslash
A_\beta| < \omega$ for all $\alpha \in a$, hence there is a fixed $s
\in [\omega]^{<\omega}$ such that $$b = \{ \alpha \in a : A_\alpha
\backslash A_\beta \subs s \}$$ is uncountable. But clearly, $s \cup
A_\beta$ is then an upper bound of $\mathcal{K} = \{ I_\alpha :
\alpha \in b \}$ in $\mathcal{I}$. Consequently, by theorem
\ref{ter}, the space $X_{\mathcal{I}}$ is as required.
\end{proof}

This result takes care of problems 4.6 and 4.7 from \cite{T} by
showing that it is consistent to have first countable Tychonov
spaces with caliber $\omega_1$ (and hence also CCC) without any
point-countable $\pi$-base. With some further elaboration we may
find examples that, in addition, are also hereditarily Lindel\"of,
and thus provide a solution to problem 4.3 as well.

\begin{theorem}\label{hl}
Let $\{ A_\alpha : \alpha < \omega_2 \} \subs \mathcal{P}(\omega)$
be a mod finite strictly increasing $\omega_2$-sequence with the
additional property that for every uncountable index set $a \in
[\omega_2]^{\omega_1}$ there is a pair $\{ \alpha,\beta \} \in
[a]^{2}$ such that $A_\alpha \subs A_\beta$, (i. e. $A_\alpha$ is
really a subset of $A_\beta$, not just mod finite). Then, for
$\mathcal{I}$ defined as in corollary \ref{inc}, the space
$X_{\mathcal{I}}$ is also hereditarily Lindel\" of.
\end{theorem}

\begin{proof}
Assume, on the contrary, that $X_{\mathcal{I}}$ has an uncountable
right-separated subspace. Without loss of generality this may be
taken of the form $\{ I_\alpha: \alpha \in a \}$, right separated in
the natural well-ordering of its indices, where $a \in
[\omega_2]^{\omega_1}$ and $I_\alpha =^* A_{\alpha}$ for all $\alpha
\in a$. Moreover, we may assume that we have a fixed $U \in
\mathcal{U}$ such that $[I_\alpha,U)_{\mathcal{I}}$ is a right
separating neighbourhood of $I_\alpha$ for any $\alpha \in a$.

Clearly, there is a fixed finite set $s \in [\omega]^{< \omega}$
such that $$b = \{ \alpha \in a : I_\alpha \Delta A_{\alpha} = s
\}$$ is uncountable. By our assumption, there  is a pair $\{\alpha,
\beta\} \in [b]^2$ (with $\alpha < \beta$) for which $A_{\alpha}
\subs A_\beta$ and hence $I_\alpha \subs I_\beta$. This, however,
would imply $I_\beta \in [I_\alpha,U)_{\mathcal{I}}$, contradicting
that $[I_\alpha,U)_{\mathcal{I}}$ is a right separating
neighbourhood of $I_\alpha$.
\end{proof}

Note that a space as in theorem \ref{hl} is a first countable
L-space, hence unlike the spaces in corollary \ref{inc}, it does not
exist under $MA_{\omega_1}$, see \cite{Sz}. Luckily, there is a
``natural" forcing construction that produces a mod finite strictly
increasing $\omega_2$-sequence in $\mathcal{P(\omega)}$ with the
additional property required in theorem \ref{hl}.

\begin{theorem}\label{forc}
There is a CCC forcing that, to any ground model, adds a mod finite
strictly increasing sequence $\{ A_\alpha : \alpha < \omega_2 \}
\subs \mathcal{P}(\omega)$ in any uncountable subsequence of which
there are two members with proper inclusion.
\end{theorem}

\begin{proof}
Let $\mathbb{P}$ consist of those finite functions $p \in
Fn(\omega_2 \times \omega, 2)$ for which $\dom(p) = a \times n$ with
$a \in [\omega_2]^{< \omega}$ and $n < \omega$. We define  $p' \le
p$ (i. e. $p'$ extends $p$) as follows: $p' \supset p$, moreover
$p'(\alpha,i) = 1$ implies $p'(\beta,i) = 1$ whenever $\alpha, \beta
\in a$ with $\alpha < \beta$ and $i \in n' \setminus n$ (of course,
here $\dom(p) = a \times n$ and $\dom(p') = a' \times n'\,$). It is
straight-forward to show that $\langle \mathbb{P}, \le \rangle$ is a
CCC notion of forcing.

Let $G \subs \mathbb{P}$ be generic, then it follows from standard
density arguments that $g = \cup G$ maps $\omega_2 \times \omega$
into $2$ and if we set $$A_\alpha = \{ i < \omega : g(\alpha,i) =
1\}$$ then $\{ A_\alpha : \alpha < \omega_2 \}$ is mod finite
strictly increasing.

To finish the proof, let us assume that $p \in \mathbb{P}$ forces
that $\dot{h}$ is an order preserving injection of $\omega_1$ into
$\omega_2$. It will clearly suffice to show that $p$ then has an
extension $q$ which forces $A_{\dot{h}(\xi)} \subs
A_{\dot{h}(\eta)}$ for some $\xi < \eta < \omega_1$.

To see this, let us choose first for each $\xi < \omega_1$ a
condition $p_\xi \le p$ and an ordinal $\alpha_\xi < \omega_2$ such
that $p_\xi \Vdash \dot{h}(\xi) = \alpha_\xi$. We may assume without
any loss of generality that for some $n < \omega$ we have
$\dom(p_\xi) = a_\xi \times n$ and $\alpha_\xi \in a_\xi$ for all
$\xi$. Using standard $\Delta$-system and counting arguments, it is
easy to find then $\xi < \eta <\omega_1$ such that $p_\xi$ and
$p_\eta$ are compatible as functions and for any $i < n$ we have
$p_\xi(\alpha_\xi,i) = p_\eta(\alpha_\eta,i)$. But then we clearly
have $q = p_\xi \cup p_\eta \in \mathbb{P}$ and $q \le p$, moreover
it is obvious that $q$ forces $A_{\alpha_\xi} \subs A_{\alpha_\eta}$
and hence $A_{\dot{h}(\xi)} \subs A_{\dot{h}(\eta)}$ as well.
\end{proof}

>From theorems \ref{hl} and \ref{forc} we immediately obtain a joint
solution to problems 4.3 and 4.7 (and hence 4.6) of Tka\v cuk from
\cite{T}.

\begin{corollary}
It is consistent to have a first countable, hereditarily Lindel\"of
0-dimensional $T_2$ space with caliber $\omega_1$ having no
point-countable $\pi$-base.
\end{corollary}

Let us recall here that the failure of CH is not sufficient to
produce a mod finite strictly increasing $\omega_2$-sequence in
$\mathcal{P(\omega)}$, the basic ingredient of our examples in this
section. In fact, Kunen had proved (see e. g. \cite{JSSz}) that if
one adds $\omega_2$ Cohen reals to a model of CH then no such
sequence exists in the extension. Actually, we can show, as a
strengthening of this, that in the same model if a subfamily
$\mathcal{I}$ of $\mathcal{P}(\omega)$ has $\omega_1$ as a set
caliber then $\cof(\mathcal{I}) \le \omega$. So, we may not hope to
use the methods of this section to find examples just assuming the
negation of CH. So the following problem can be raised.

\begin{problem}
Does $2^\omega > \omega_1$ imply the existence of a first countable
Lindel\"of and/or CCC Tychonov space having no point-countable
$\pi$-base?
\end{problem}

\end{document}